\documentclass[]{eptcs}

\input{header}

\title{DisCoPy: Monoidal Categories in Python}
\author{Giovanni de Felice, Alexis Toumi, Bob Coecke
\institute{Department of Computer Science, University of Oxford.
\quad Cambridge Quantum Computing Ltd.}
\email{\{firstname.lastname\}@cs.ox.ac.uk}}

\def\titlerunning{DisCoPy: Monoidal Categories in Python}
\def\authorrunning{G. de Felice, A. Toumi \& B. Coecke}

\begin{document}
\maketitle


\begin{abstract}
We introduce DisCoPy, an open source toolbox for computing with monoidal categories.
The library provides an intuitive syntax for defining string diagrams and monoidal functors.
Its modularity allows the efficient implementation of computational experiments in the various applications of category theory where diagrams have become a lingua franca.
As an example, we used DisCoPy to perform natural language processing on quantum hardware for the first time.
\end{abstract}

\section*{Introduction} 

String diagrams are a graphical calculus for monoidal categories, introduced independently by Hotz \cite{Hotz65} in computer science and Penrose \cite{Penrose71} in physics, then formalised by Joyal and Street \cite{JoyalStreet88, JoyalStreet91}.
Graphical languages are surveyed in Selinger \cite{Selinger10},
they have become a standard tool in applied category theory.
We cite a few of the growing list of applications:
computer science \cite{BrownHutton94, Abramsky96},
quantum theory and the ZX-calculus \cite{Coecke05,AbramskyCoecke08,CoeckeDuncan08},
networks and control theory \cite{BaezErbele14, BaezFong15, BaezPollard17}, concurrency \cite{BonchiEtAl14a},
databases and
knowledge representation \cite{Patterson17, BonchiEtAl18},
Bayesian reasoning and causality \cite{CoeckeSpekkens12,ChoJacobs19,KissingerUijlen19},
linguistics and cognition \cite{ClarkEtAl08,BoltEtAl17},
functional programming \cite{Riley18},
machine learning and
game theory \cite{FongEtAl17, GhaniEtAl18}.
In all these applications, string diagrams are the syntax and structure-preserving functors are used to compute their semantics in concrete categories.

There are several existing proof assistants with graphical user interfaces for rewriting string diagrams in a more or less automated fashion:
quantomatic \cite{KissingerZamdzhiev15} and PyZX \cite{KissingervandeWetering19} for the ZX-calculus,
globular \cite{BarEtAl} and its successor homotopy.io \cite{ReutterVicary19} for higher categories,
cartographer \cite{SobocinskiEtAl19} for symmetric monoidal categories.
However, these are all stand-alone tools which use different task-specific encodings for diagrams, preventing interoperability between them and with the software ecosystems of application domains.

DisCoPy (Distributional Compositional Python) is not yet-another rewriting tool, rather it is meant as a toolbox for compiling diagrams into code, be it for neural networks, tensor computation or quantum circuits.
It provides an intuitive Python syntax for diagrams, allowing to visualise and reason about computation at a high level of abstraction.
Monoidal functors allow to translate these diagrams into concrete computation, interfacing with optimised task-specific libraries.
DisCoPy is an open source package, it is available with an extensive documentation and demonstration notebooks hosted at:

\url{https://github.com/oxford-quantum-group/discopy}

This paper describes the architecture of DisCoPy, focusing on the translation from abstract categorical definitions into their concrete implementation in Python.
We assume some working knowledge of category theory and refer the reader to \cite{Lane98} and to \cite{Awodey06} for an introduction.
Implementing a category in an object-oriented programming language amounts to defining a pair of classes for its objects and arrows, together with a pair of methods for identity and composition.
When the category is free, composition is implemented by list concatenation and identity by the empty list.
Concrete categories may then be defined by subclassing this free category and overriding identity and composition.
These are expected to respect the usual associativity and unit axioms, however they cannot be formally checked in Python.

Starting from free categories (section~\ref{1-cat}) as a base class, more structure can be added by subclassing and adding new methods.
Quotient categories can be implemented by a method for computing normal forms.
For instance, monoidal categories (section~\ref{2-monoidal}) subclass categories with an extra method for tensor product and one for interchanger normal form.
For now, we implemented Cartesian and rigid monoidal categories (section~\ref{4-rigid}), as these provide a syntax for the concrete categories implemented in DisCoPy: Python functions (appendix~\ref{3-cartesian}) and numpy \cite{VanderWaltEtAl11} tensors (section~\ref{5-tensor}).
The development of DisCoPy was first motivated by the implementation of natural language processing on quantum hardware.
Hence, we implemented quantum circuits (section~\ref{6-circuit}) as a subclass of rigid monoidal categories with an extra method for evaluation as numpy tensors and interface with the t|ket$\rangle$ compiler \cite{SivarajahEtAl20}.

We hope that this toolbox will prove to be of use to the applied category theory community, and plan on adding more categorical tools to it in the future.

\section{cat.py}\label{1-cat}

We give a brief introduction to the concepts of signatures, free categories
and functors, fixing some notation for the rest of the paper.
We then describe how these concepts, once implemented in Python, give the
architecture of \py{cat}, the first module of the DisCoPy package.

A \emph{simple signature} $\Sigma$ is given by a pair of sets $\Sigma_0, \Sigma_1$ of generating objects and arrows, and a pair of functions
$\mathtt{dom}, \mathtt{cod} : \Sigma_1 \to \Sigma_0$ called the domain and the codomain
respectively.
A morphism of signatures $F : \Sigma \to \Sigma'$ is a pair of functions
$F_0 : \Sigma_0 \to \Sigma'_0$, $F_1 : \Sigma_1 \to \Sigma'_1$ which commute with $\mathtt{dom}$ and $\mathtt{cod}$.
The \emph{free category} $\mathbf{C}(\Sigma)$ generated by a signature $\Sigma$, is
defined as follows: the objects are given by $\Sigma_0$ and an arrow
$f : s \to t$ is given by a list $f = f_1 ... f_n \in \Sigma_1^n$ such that
$\mathtt{dom}(f_1) = s, \mathtt{cod}(f_n) = t$ and $\mathtt{cod}(f_i) = \mathtt{dom}(f_{i + 1})$
for all $i \leq n$. Identity arrows are given by the empty list and composition
is given by list contatenation.
The universal property of $\mathbf{C}(\Sigma)$ may be stated as follows: functors
$F : \mathbf{C}(\Sigma) \to \mathbf{D}$ from a free category are uniquely defined by
their image on the signature $\Sigma$, i.e. by a morphism of signatures
$F : \Sigma \to U(\mathbf{D})$ for $U$ the forgetful functor.
These abstract definitions are implemented in an object-oriented fashion with
the following Python classes, where the tuples denote the arguments to the
corresponding class constructors.

\begin{class}\normalfont\texttt{cat.Ob(name)} is given by any Python object \py{name}.
Equality of objects is given by the equality of names, i.e. \py{x == y}
if and only if \py{x.name == y.name}.
\end{class}

\begin{class}\normalfont\texttt{cat.Arrow(dom, cod, boxes)} is given by a pair of \py{Ob}-instances
\py{dom, cod} and a sequence of \py{Box}-instances \py{boxes}.
Axioms are checked at initialisation, \py{AxiomError} is raised otherwise.
An \py{Arrow}-instance \py{f} has the following methods:

\begin{itemize}
\item \py{f.dagger()} returns the dagger of an arrow,
it has syntactic sugar \py{f[::-1]}.
\item \py{f.then(g)} returns the composition of two arrows,
it is written \py{f >> g} or \py{g << f}.
\item \py{Arrow.id(x)} returns the identity arrow
on a given object, it is shortened to \py{Id(x)}.
\end{itemize}
Axioms are implemented as part of the testing suite for the DisCoPy
package: associativity \py{f >> (g >> h)}
\s \py{== (f >> g) >> h} and unitality
\py{f >> Id(f.cod) == f == Id(f.dom) >> f} follow directly from the monoid
structure of Python sequences.
\py{Arrow}-instances can be manipulated using the standard Python syntax for
sequences: they are iterable, indexable and sliceable, i.e.
\py{f[i] == f.boxes[i]} and \py{f[i:j]} returns the arrow \py{f[i] >> ...}
\py{>> f[j]}.
Printing an arrow yields a DisCoPy expression which generates it,
which proves useful for debugging and interactive programming.
\end{class}

\begin{class}\normalfont\texttt{cat.Box(name, dom, cod, data=None, \_dagger=False)} is a subclass of \py{Arrow}, i.e. a box \py{f} is defined as an arrow with \py{f.boxes == [f]}.
\py{data} is an optional argument, it can be used to attach arbitrary Python data to a box.
\py{\_dagger} is an optional Boolean argument, set to \py{False} by default.
It is used to construct an involutive identity-on-objects contravariant endo-functor, i.e. a dagger.
\end{class}

\begin{class}\normalfont\texttt{cat.Functor(ob, ar, ob\_factory=Ob, ar\_factory=Arrow)} is given by a pair of mappings \py{ob, ar} from \py{Ob} to \py{ob\_factory} and from \py{Box} to \py{ar\_factory} respectively.
The domain of a \py{Functor}-instance \py{F} is defined implicitly as the free category generated by the domain of \py{ob} and \py{ar}.
The two factory methods allow to build functors with arbitrary codomains, \py{ar\_factory} is required to provide methods \py{id} and \py{then}.
\end{class}

\begin{example}
In order to illustrate the syntax, we give a basic example of a free functor.
\normalfont
\begin{python}
x, y, z = Ob('x'), Ob('y'), Ob('z')
f, g, h = Box('f', x, y), Box('g', y, z), Box('h', z, x)
F = cat.Functor(ob={x: y, y: z, z: x}, ar={f: g, g: h})
assert F(f >> g) == F(f) >> F(g) == g >> h
\end{python}
\end{example}

\begin{example}
\py{PythonFunctor} implements functors into Python functions, see
appendix~\ref{3-cartesian}.
\normalfont
\begin{python}
F = PythonFunctor(
    ob={x: 0, y: 1, z: 2}, ar={f: lambda: 42, g: lambda x: (x, x + 1)})
assert F(f >> g)() == (42, 43)
\end{python}
\end{example}

\begin{example}
The category of matrices is implemented by the \py{Tensor} class, see
section~\ref{5-tensor}.
\normalfont
\begin{python}
F = TensorFunctor(ob={x: 1, y: 2, z: 2}, ar={f: [0, 1], g: [0, 1, 1, 0]})
assert F(f >> g) == F(f) >> F(g) == [1, 0]
\end{python}
\end{example}

All the objects, arrows, boxes and functors constructed with DisCoPy are representable in the following sense: for any instance \py{x}, the string \py{repr(x)} is a valid Python expression which evaluates to \py{x}, assuming that the names defining \py{x} are representable themselves.
Hence, DisCoPy data structures can be serialized and exported in a standard format \cite{Statebox20} for interoperability with other category theory tools.
Note that while the data structures and methods of DisCoPy are purely functional, boxes may hold mutable a \py{data} attribute.
It is used to build boxes that are not finitely generated such as the weights of neural networks (appendix~\ref{3-cartesian}) or the phases of quantum circuits (section~\ref{6-circuit}).

\section{monoidal.py}\label{2-monoidal}

In this section, we describe the core data structure behind DisCoPy:
the \py{Diagram} class, an implementation of the arrows of a free
monoidal category. We begin with a definition of free monoidal categories via
free \emph{premonoidal categories}.

\begin{definition}\label{definition-1}
A (strict) monoidal category is a category $\mathbf{C}$ equipped with an associative
and unital functor $\otimes : \mathbf{C} \times \mathbf{C} \to \mathbf{C}$ for $\times$
the Cartesian product.
A (strict) premonoidal category is a category $\mathbf{C}$ equipped with an
associative and unital functor $\boxtimes : \mathbf{C} \ \Box \ \mathbf{C} \to \mathbf{C}$
for $\Box$ given by the following pushout in $\mathbf{Cat}$:

\begin{center}\begin{tikzcd}
& \mathbf{C}_0 \times \mathbf{D}_0 \arrow[r, ""] \arrow[d, ""]
&  \mathbf{C} \times \mathbf{D}_0 \arrow[d, ""] \\
&  \mathbf{C}_0 \times \mathbf{D} \arrow[r, ""]
&  \mathbf{C} \ \Box \  \mathbf{D}
\end{tikzcd}\end{center}
where $\mathbf{C}_0, \mathbf{D}_0$ are the discrete category of objects and the maps are given by the inclusions.
\end{definition}

\begin{example}\label{example-2}
The Keisli category for a strong monad over a monoidal category is a
premonoidal category. It is monoidal precisely when the monad is commutative.
The state monad over the category of sets yields a denotational
semantics for side-effects \cite{PowerRobinson97}.
\end{example}

\begin{example}
For a semiring $\bb{S}$, the distribution monad $X \mapsto \bb{S}^X$ over the
category of finite sets yields a premonoidal category with matrices over $\bb{S}$ as arrows and Kronecker product as tensor. It is monoidal when $\bb{S}$ is commutative.
\end{example}

The pushout $\Box$ is known as the \emph{funny tensor product}, it may also be defined
by the following universal property. Let $\mathbf{C} \Rightarrow \mathbf{D}$ be the
category where objects are functors $\mathbf{C} \to \mathbf{D}$ and arrows are the
transformations with no naturality requirement. Then $-\ \Box\ \mathbf{C}$
is characterised as the left adjoint of $\mathbf{C} \Rightarrow -$ in $\mathbf{Cat}$.
This makes $\Box$ a closed symmetric monoidal structure over $\mathbf{Cat}$, the
unique such structure apart from the usual Cartesian product, see
\cite{FoltzEtAl80}. A functor from the funny tensor product
$\mathbf{C}\ \Box\ \mathbf{D}$ can be understood as a functor which is ``separately
functorial'' in its two arguments $\mathbf{C}$ and $\mathbf{D}$, in analogy to separate continuity.
Premonoidal categories are also known as one-object sesquicategories
(one-and-a-half categories), i.e. 2-categories without the interchange law,
see \cite{Street96}.
Our motivation for working in a premonoidal setting is two-fold: 1) free
premonoidal categories have a simple presentation as free categories and 2)
free monoidal categories may then be described as quotient categories.

We define a \emph{monoidal signature} $\Sigma$ as a pair of sets $\Sigma_0, \Sigma_1$ with a pair of functions into the free monoid $\mathtt{dom}, \mathtt{cod} : \Sigma_1 \to \Sigma_0^\star$.
Generating arrows $f \in \Sigma_1$ are depicted as boxes with input $\mathtt{dom}(f)$ and output $\mathtt{cod}(f)$.
Given a monoidal signature $\Sigma$, we define a simple signature $L(\Sigma)$ with objects $\Sigma_0^\star$, arrows $\Sigma_0^\star \times \Sigma_1 \times \Sigma_0^\star$ with $\mathtt{dom}(u, f, v) = usv$
and $\mathtt{cod}(u, f, v) = utv$ for $s = \mathtt{dom}(f)$ and $t = \mathtt{cod}(f)$.
A layer $(u, f, v) \in \Sigma_0^\star \times \Sigma_1 \times \Sigma_0^\star$ is depicted as a box with wires to its left and right:

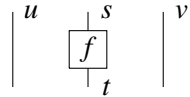
\begin{center}
\begin{tikzpicture}[baseline=(O.base)]
\node (O) at (0, 0.5) {};
\node () at (0.25, 1) {$u$};
\node () at (1.25, 1) {$s$};
\node () at (2.25, 1) {$v$};
\node () at (1.25, 0.0) {$t$};
\draw [out=-90, in=90] (0, 1) to (0, 0);
\draw [out=-90, in=90] (1, 1) to (1, 0.75);
\draw [out=-90, in=90] (2, 1) to (2, 0);
\draw [out=-90, in=90] (1.0, 0.25) to (1.0, 0);
\draw (0.75, 0.25) -- (1.25, 0.25) -- (1.25, 0.75) -- (0.75, 0.75) -- (0.75, 0.25);
\node () at (1.0, 0.5) {$f$};
\end{tikzpicture}
\end{center}

\begin{proposition}
Given a monoidal signature $\Sigma$, the free premonoidal category is the free category $\mathbf{PMC}(\Sigma) = \mathbf{C}(L(\Sigma))$ generated by the simple signature of layers $L(\Sigma)$.
\end{proposition}

We define a \emph{diagram} as an arrow of $\mathbf{C}(L(\Sigma))$.
Diagrams are uniquely defined by a domain, a list of generators and an \emph{offset} for each box: the number of wires passing to its left.
Two diagrams are equal in the free monoidal category if they are related by a series of \emph{interchangers},
where $u, v, w \in \Sigma_0^\star$ and $s \xto{f} t, s' \xto{f'} t' \in \Sigma_1$:

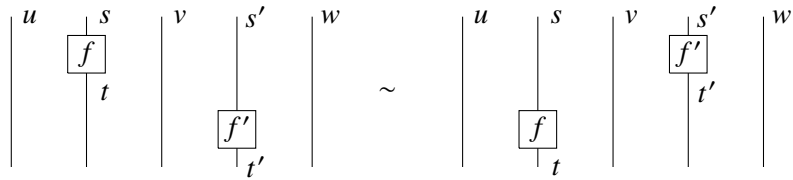
\begin{center}
\begin{tikzpicture}[baseline=(O.base)]
\node (O) at (0, 1.0) {};
\node () at (0.25, 2.0) {$u$};
\node () at (1.25, 2.0) {$s$};
\node () at (2.25, 2.0) {$v$};
\node () at (3.25, 2.0) {$s'$};
\node () at (4.25, 2.0) {$w$};
\node () at (1.25, 1.0) {$t$};
\node () at (3.25, 0.0) {$t'$};
\draw [out=-90, in=90] (0, 2.0) to (0, 0.0);
\draw [out=-90, in=90] (1, 2.0) to (1, 1.75);
\draw [out=-90, in=90] (2, 2.0) to (2, 0.0);
\draw [out=-90, in=90] (3, 2.0) to (3, 0.75);
\draw [out=-90, in=90] (4, 2.0) to (4, 0.0);
\draw [out=-90, in=90] (1.0, 1.25) to (1.0, 0.0);
\draw [out=-90, in=90] (3.0, 0.25) to (3.0, 0.0);
\draw (0.75, 1.25) -- (1.25, 1.25) -- (1.25, 1.75) -- (0.75, 1.75) -- (0.75, 1.25);
\node () at (1.0, 1.5) {$f$};
\draw (2.75, 0.25) -- (3.25, 0.25) -- (3.25, 0.75) -- (2.75, 0.75) -- (2.75, 0.25);
\node () at (3.0, 0.5) {$f'$};
\node () at (5, 1.0) {$\sim$};
\node () at (6.25, 2.0) {$u$};
\node () at (7.25, 2.0) {$s$};
\node () at (8.25, 2.0) {$v$};
\node () at (9.25, 2.0) {$s'$};
\node () at (10.25, 2.0) {$w$};
\node () at (9.25, 1.0) {$t'$};
\node () at (7.25, 0.0) {$t$};
\draw [out=-90, in=90] (6, 2.0) to (6, 0.0);
\draw [out=-90, in=90] (7, 2.0) to (7, 0.75);
\draw [out=-90, in=90] (8, 2.0) to (8, 0.0);
\draw [out=-90, in=90] (9, 2.0) to (9, 1.75);
\draw [out=-90, in=90] (10, 2.0) to (10, 0.0);
\draw [out=-90, in=90] (9.0, 1.25) to (9.0, 0.0);
\draw [out=-90, in=90] (7.0, 0.25) to (7.0, 0.0);
\draw (8.75, 1.25) -- (9.25, 1.25) -- (9.25, 1.75) -- (8.75, 1.75) -- (8.75, 1.25);
\node () at (9.0, 1.5) {$f'$};
\draw (6.75, 0.25) -- (7.25, 0.25) -- (7.25, 0.75) -- (6.75, 0.75) -- (6.75, 0.25);
\node () at (7.0, 0.5) {$f$};
\end{tikzpicture}
\end{center}

The \emph{right interchangers}, going from the right- to the left-hand side of the previous equivalence are terminating on boundary-connected diagrams.
Given a boundary-connected diagram with $n$ boxes, a normal form can be reached in at most $O(n^3)$ steps \cite[Theorem~36]{DelpeuchVicary18}.
This makes the word problem for monoidal categories --- i.e. given two arbitrary diagrams, are they equal up to interchanger? --- decidable in polynomial time, see \cite[Theorem~48]{DelpeuchVicary18}.

\begin{proposition}
Given a monoidal signature $\Sigma$, the free monoidal category is the quotient $\mathbf{MC}(\Sigma) = \mathbf{PMC}(\Sigma) / \cal{I}$ for $\cal{I}$ the interchanger relation, see \cite[Theorem~16]{DelpeuchVicary18}.
\end{proposition}

\begin{example}
Context-free grammars (CFGs) are a special case of free monoidal categories, where
non-terminal symbols are generating objects and production rules are generators
with an atomic type as domain. Syntax trees are diagrams, their normal form
correspond to the leftmost derivation in a CFG. Weighted CFGs can be defined
as free monoidal categories equipped with a monoidal functor into a monoid delooping, see \cite{ShieblerEtAl20}.
\end{example}

Note that a monoidal category generated by one object, i.e. with
$\Sigma_0 = \set{1}$, is called a \emph{PRO}. An example is the category of
circuits described in section~\ref{6-circuit}.
This combinatorial definition yields an implementation of the free monoidal category where the generating objects $\Sigma_0$ are given by \py{Ob}-instances and the generating arrows $\Sigma_1$ are given by \py{Box}-instances.

\begin{class}\normalfont\texttt{monoidal.Ty(x\_1, ..., x\_n)} is a subclass of \py{cat.Ob}.
A type \py{x} is given by a (possibly empty) list of objects,
It has a method \py{x.tensor(y)}, shortened to \py{x @ y}, which inherits the monoid structure of lists with the empty type \py{Ty()} as unit.
\end{class}

\begin{class}\normalfont\texttt{monoidal.PRO(n)} is a subclass of \py{Ty} generated by \py{Ob(1)}.
A \py{PRO} type is given a natural number \py{n} with addition as tensor, i.e.
\py{PRO(n) == Ty(n * (1, ))}.
\end{class}

\begin{class}\normalfont\texttt{monoidal.Layer(left, box, right)} is a subclass of \py{cat.Box}.
A layer is given by a \py{Box}-instance \py{box} and a pair of \py{Ty}-instances \py{left} and \py{right}.
While essential to the internal structure of DisCoPy, layers remain invisible to the end user of the package.
\end{class}

\begin{class}\normalfont\texttt{monoidal.Diagram(dom, cod, boxes, offsets, layers=None)} subclasses \py{Arrow}.
A diagram \py{f} is given by a pair of types \py{dom, cod}, a pair of equal-length sequences \py{boxes, offsets} of \py{Box}-instances and natural numbers respectively.
\py{layers} is an optional argument, if omitted it will be computed from \py{dom}, \py{boxes} and \py{offsets}.
\py{f} has all the methods of \py{Arrow} plus the following:
\begin{itemize}
\item \py{f.tensor(g)}, shortened to \py{f @ g}, returns the diagram with
\\\py{(f @ g).boxes == f.boxes + g.boxes} and
\\\py{(f @ g).offsets == f.offsets + [n + len(f.cod) for n in g.offsets]}.
\item \py{f.interchange(i, j, left=False)} returns a diagram with the boxes at indices \py{i} and \py{j} interchanged or raises an \py{InterchangerError} if they are connected.
This method gets called recursively whenever \py{i < j + 1} or \py{j < i - 1}.
If there is a choice in how to interchange, then we apply the right interchanger by default.
\item \py{f.normalize(left=False)} yields a reduction sequence applying right interchangers repeatedly. Setting \py{left=True} will apply left interchangers instead.
\item \py{f.normal\_form(left=False)} returns the last output of \py{f.normalize} if \py{f} is boundary-connected, otherwise it raises a \py{NotImplementedError}.
\item \py{f.draw(to\_tikz=False, **params)} draws a diagram using \texttt{networkx} \cite{Networkx20}, see the DisCoPy documentation \cite{DeFeliceToumi} for a complete list of parameters. Diagrams are displayed with \texttt{matplotlib} \cite{Matplotlib20} by default, setting \py{to\_tikz=True} outputs a list of Ti\textit{k}z \cite{Tantau13} commands instead. This is how all the diagrams in this article were produced.
\item \py{f.to\_gif(g, h, ..., **params)} takes a list of diagrams and returns an animated GIF which can be used to visualise a rewriting process in a jupyter notebook \cite{KluyverEtAl16}.
\end{itemize}

As for \py{Arrow}, \py{Diagram}-instances can be manipulated as lists.
Indexing a diagram \py{f} returns a layer \py{Id(left) @ box @ Id(right)} such that \py{f == f[0] >> ... >> f[-1]}.
\end{class}

\begin{class}\normalfont\texttt{monoidal.Box(name, dom, cod, data=None, \_dagger=False)} is a subclass of \py{cat.Box} and \py{Diagram}. A box \py{f} is defined as a diagram with \py{f.boxes, f.offsets == [f], [0]}.
\end{class}

\begin{class}\normalfont\texttt{monoidal.Functor(ob, ar, ob\_factory=Ty, ar\_factory=Diagram)} is a subclass of \py{cat.} \py{Functor}.
A monoidal functor \py{F} is given by a pair of mappings \py{ob, ar} from \py{Ty} to \py{ob\_factory} and from \py{Box} to \py{ar\_factory} respectively.
Factory methods are optional arguments that allow arbitrary codomains.
\py{ob\_factory} is required to provide a \py{tensor} method, \py{ar\_factory} should provide \py{id}, \py{then} and \py{tensor}.
\end{class}

The implementation of the \py{Diagram} data structure and its normal form follows directly from their formal definitions, they require no further explanation.
The implementation of the \py{draw} method however requires some non-trivial choices: given the combinatorial encoding of a diagram, which embedding on the plane should we return?
Moreover, the drawing algorithm may be treated as a proof of the equivalence from the combinatorial to the geometric definition of diagrams introduced by Joyal and Street \cite{JoyalStreet88, JoyalStreet91}.
The other direction, i.e. from a planar embedding (encoded as a grid of pixels) to its combinatorial encoding as a DisCoPy \py{Diagram}, is part of an in-development application for automated string diagram recognition.
We leave the details to appendix~\ref{a-draw}.

\section{rigid.py}\label{4-rigid}

In this section, we describe the module \py{rigid}, an implementation of
rigid monoidal categories.
We present \emph{snake removal}, the algorithm for normalising rigid diagrams, and its application to the semantics of \emph{pregroup grammars}.
Note that rigid categories are also called \emph{autonomous} \cite{JoyalStreet88,Delpeuch14a}, they are equivalent to the compact 2-categories of Preller, Lambek \cite{PrellerLambek07} with one object.

A (strict) monoidal category $\mathbf{C}$ is \emph{rigid} when every object $x$
has left and right adjoints $x^l$ and $x^r$ and four morphisms
$x \otimes x^l \xto{\epsilon} 1 \xto{\eta} x^l \otimes x$ and
$x^r \otimes x \xto{\epsilon'} 1 \xto{\eta'} x \otimes x^r$ depicted as cups and caps,
subject to $(\epsilon' \otimes 1_x) \circ (1_x \otimes \eta') = 1_x = (1_x \otimes \epsilon) \circ (\eta \otimes 1_x)$.
Note that any monoidal functor between rigid categories is isomorphic to a rigid functor, i.e. which sends cups to cups and caps to caps.

In a rigid category, left and right adjoints are unique up to a unique isomorphism.
They cancel each other --- i.e. $(x^l)^r = x = (x^r)^l$ --- and they are anti-homomorphisms of the tensor product --- i.e. $1^l = 1^r = 1$, $(x \otimes y)^l = y^l \otimes x^l$ and $(x \otimes y)^r = y^r \otimes x^r$.
Thus, the adjoint of a product can always be written as a product of adjoints.
Given a set of atomic types $\Sigma_0$, the objects of the free rigid category are given by lists of pairs $(x, n) \in \Sigma_0 \times \bb{Z}$ with the inclusion $\Sigma_0 \injects \Sigma_0 \times \bb{Z}$ defined by $x \mapsto (x, 0)$ and the adjoints $(x, n)^l = (x, n - 1)$ and $(x, n)^r = (x, n + 1)$.
We define a \emph{rigid signature} $\Sigma$ as a pair of sets $\Sigma_0, \Sigma_1$ with a pair of functions $\mathtt{dom}, \mathtt{cod} : \Sigma_1 \to (\Sigma_0 \times \bb{Z})^\star$,
then the free rigid category is given by the quotient $\mathbf{RC}(\Sigma) = \mathbf{MC}(\Sigma') / \cal{R}$
where $\Sigma'$ is the monoidal signature with $\Sigma'_0 = \Sigma_0 \times \bb{Z}$ and $\Sigma'_1 = \Sigma_1 + \set{\mathtt{cup}_x : x \otimes x^r \to 1, \mathtt{cap}_x : 1 \to x \otimes x^l}_{x \in \Sigma'_0}$.
The cups and caps for the unit are the identity, those for $x \otimes y$ are given by nesting the cups and caps of $x$ and $y$.
The relations $\cal{R}$ are given by the snake equation for every $x \in \Sigma'_0$:
\begin{center}
\begin{tikzpicture}[baseline=(O.base)]
\node (O) at (0, 1.0) {};
\node () at (0.25, 2.0) {$x$};
\node () at (1.25, 1.0) {$x^r$};
\node () at (2.25, 1.0) {$x$};
\draw [out=-90, in=90] (0, 2.0) to (0, 0.75);
\draw [out=180, in=90] (1.5, 1.5) to (1.0, 1.25);
\draw [out=0, in=90] (1.5, 1.5) to (2.0, 1.25);
\draw [out=-90, in=90] (1.0, 1.25) to (1.0, 0.75);
\draw [out=-90, in=90] (2.0, 1.25) to (2.0, 0.0);
\draw [out=-90, in=180] (0, 0.75) to (0.5, 0.5);
\draw [out=-90, in=0] (1.0, 0.75) to (0.5, 0.5);
\node () at (3.0, 1.0) {$\sim$};
\node () at (4.25, 2.0) {$x$};
\draw [out=-90, in=90] (4.0, 2.0) to (4.0, 0.0);
\node () at (5.0, 1.0) {$\sim$};
\node () at (8.25, 2.0) {$x$};
\node () at (6.25, 1.0) {$x$};
\node () at (7.25, 1.0) {$x^l$};
\draw [out=-90, in=90] (8.0, 2.0) to (8.0, 0.75);
\draw [out=180, in=90] (6.5, 1.5) to (6.0, 1.25);
\draw [out=0, in=90] (6.5, 1.5) to (7.0, 1.25);
\draw [out=-90, in=90] (6.0, 1.25) to (6.0, 0.0);
\draw [out=-90, in=90] (7.0, 1.25) to (7.0, 0.75);
\draw [out=-90, in=180] (7.0, 0.75) to (7.5, 0.5);
\draw [out=-90, in=0] (8.0, 0.75) to (7.5, 0.5);
\end{tikzpicture}
\end{center}

The \py{rigid} module implements classes \py{Ty}, \py{Diagram}, \py{Box} and \py{Functor} which subclass those from \py{monoidal}.
Pregroup types, i.e. the objects of free rigid categories, are implemented as \py{Ty}-instances \py{x} with two attributes \py{x.l} and \py{x.r} for the adjoints.
\py{Box} has two subclasses \py{Cup} and \py{Cap} implementing the adjunction for simple types.
The \py{Diagram} class has two static methods \py{cups} and \py{caps} which implement the adjunction for product types.
The \py{normalize} and \py{normal\_form} methods are overriden to implement snake removal: for each pair of cup and cap forming a snake, we first apply interchangers to make them adjacent, then replace the snake with an identity, see \cite[Definition 2.12]{DunnVicary19}.

\begin{example}
A rigid category is compact-closed if it is also symmetric monoidal.
In that case, the left and right adjoints coincide.
The category of matrices over a commutative semiring with Kronecker product as tensor is compact-closed, thus it is rigid.
\end{example}

\begin{example}\label{example-pregroup}
Lambek's pregroup grammars \cite{Lambek99,Lambek01,Lambek08} can be defined in terms of free rigid categories.
Indeed, a pregroup is a thin rigid category, i.e. with at most one arrow between any two objects.
A pregroup grammar $G$ is given by a vocabulary $V$, a finite set of basic types $B$ with $s \in B$ the sentence type and a finite dictionnary $D \sub V \times (B \times \bb{Z})^\star$ assigning pregroup types to words.
Let $\Sigma_G$ be the rigid signature with generating objects $V + B$ and arrows $w \to t$ for each dictionnary entry $(w, t) \in D$ and $\mathbf{G} = \mathbf{RC}(\Sigma_G)$.
The language of $G$ is given by $L(G) = \set{u \in V^\star \ \vert \ \exists \ f : u \to s \in \mathbf{G}}$.
That is, a list of words $u \in V^\star$ is a grammatical sentence whenever there is a diagram $f = g \circ d$ with $d : u \to t$ a product of dictionnary entries and $g : t \to s$ a pregroup derivation generated by cups and caps.
We do not draw the wires for words and depict the dictionnary entries as triangles. For example, let $B = \set{s, n}$ and $D = \set{\text{one, two, three} \to n, \s \text{plus} \to n^r \otimes n \otimes n^l, \s \text{equals} \to n^r \otimes s \otimes n^l}$, then ``one plus two equals three'' is a grammatical sentence:
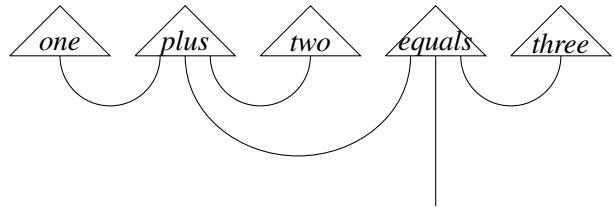
\begin{center}
\begin{tikzpicture}[baseline=(O.base), scale=0.666]
\node (O) at (0, 0) {};
\draw (0.0, 0) -- (2.0, 0) -- (1.0, 1) -- (0.0, 0);
\node () at (1.0, 0.25) {one};
\draw (2.5, 0) -- (4.5, 0) -- (3.5, 1) -- (2.5, 0);
\node () at (3.5, 0.25) {plus};
\draw (5.0, 0) -- (7.0, 0) -- (6.0, 1) -- (5.0, 0);
\node () at (6.0, 0.25) {two};
\draw (7.5, 0) -- (9.5, 0) -- (8.5, 1) -- (7.5, 0);
\node () at (8.5, 0.25) {equals};
\draw (10.0, 0) -- (12.0, 0) -- (11.0, 1) -- (10.0, 0);
\node () at (11.0, 0.25) {three};
\draw [out=-90, in=180] (1.0, 0) to (2.0, -1);
\draw [out=-90, in=0] (3.0, 0) to (2.0, -1);
\draw [out=-90, in=180] (4.0, 0) to (5.0, -1);
\draw [out=-90, in=0] (6.0, 0) to (5.0, -1);
\draw [out=-90, in=180] (9.0, 0) to (10.0, -1);
\draw [out=-90, in=0] (11.0, 0) to (10.0, -1);
\draw [out=-90, in=180] (3.5, 0) to (5.75, -2);
\draw [out=-90, in=0] (8.0, 0) to (5.75, -2);
\draw [out=-90, in=90] (8.5, 0) to (8.5, -3);
\end{tikzpicture}
\end{center}
\end{example}

\begin{example}
For any monoidal category $\mathbf{C}$, there is a free rigid category $A(\mathbf{C})$ with a fully-faithful monoidal functor $\mathbf{C} \injects A(\mathbf{C})$, see \cite{Delpeuch14a}.
Concretely, this means that in a rigid diagram with boxes coming from a monoidal category, if the domain and codomain have no adjoint types then all snakes can be removed.
This allows to give a semantics to a pregroup grammar $G$ as a rigid functor $F : \mathbf{G} \to A(\mathbf{C})$.
For example, let $\text{one, two} : 1 \to n$ and $\text{plus} : n \times n \to n$ be functions, then we can compute the meaning of ``one plus two'':
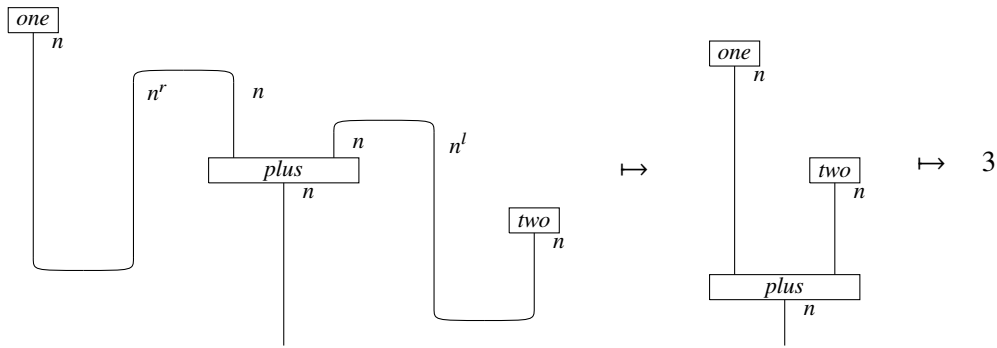
\begin{center}
\begin{tikzpicture}[baseline=(O.base), scale=0.666]
\node (O) at (0, 3.5) {};
\node [scale=0.8] () at (0.5, 6.05) {$n$};
\node [scale=0.8] () at (2.5, 5.05) {$n^r$};
\node [scale=0.8] () at (4.5, 5.05) {$n$};
\node [scale=0.8] () at (6.5, 4.05) {$n$};
\node [scale=0.8] () at (8.5, 4.05) {$n^l$};
\node [scale=0.8] () at (5.5, 3.05) {$n$};
\node [scale=0.8] () at (10.5, 2.05) {$n$};
\draw (0.0, 6.25) .. controls (0.0, 1.75) .. (0.0, 1.75);
\draw (3.0, 5.5) .. controls (2.0, 5.5) .. (2.0, 5.25);
\draw (3.0, 5.5) .. controls (4.0, 5.5) .. (4.0, 5.25);
\draw (2.0, 5.25) .. controls (2.0, 1.75) .. (2.0, 1.75);
\draw (4.0, 5.25) .. controls (4.0, 3.75) .. (4.0, 3.75);
\draw (7.0, 4.5) .. controls (6.0, 4.5) .. (6.0, 4.25);
\draw (7.0, 4.5) .. controls (8.0, 4.5) .. (8.0, 4.25);
\draw (6.0, 4.25) .. controls (6.0, 3.75) .. (6.0, 3.75);
\draw (8.0, 4.25) .. controls (8.0, 0.75) .. (8.0, 0.75);
\draw (5.0, 3.25) .. controls (5.0, 0.0) .. (5.0, 0.0);
\draw (10.0, 2.25) .. controls (10.0, 0.75) .. (10.0, 0.75);
\draw (0.0, 1.75) .. controls (0.0, 1.5) .. (1.0, 1.5);
\draw (2.0, 1.75) .. controls (2.0, 1.5) .. (1.0, 1.5);
\draw (8.0, 0.75) .. controls (8.0, 0.5) .. (9.0, 0.5);
\draw (10.0, 0.75) .. controls (10.0, 0.5) .. (9.0, 0.5);
\draw (-0.5, 6.25) -- (0.5, 6.25) -- (0.5, 6.75) -- (-0.5, 6.75) -- (-0.5, 6.25);
\node [scale=0.8] () at (0, 6.5) {one};
\draw (3.5, 3.25) -- (6.5, 3.25) -- (6.5, 3.75) -- (3.5, 3.75) -- (3.5, 3.25);
\node [scale=0.8] () at (5.0, 3.5) {plus};
\draw (9.5, 2.25) -- (10.5, 2.25) -- (10.5, 2.75) -- (9.5, 2.75) -- (9.5, 2.25);
\node [scale=0.8] () at (10.0, 2.5) {two};
\node () at (12.0, 3.5) {$\mapsto$};
\node [scale=0.8] () at (14.5, 5.383333333333333) {$n$};
\node [scale=0.8] () at (16.5, 3.0500000000000003) {$n$};
\node [scale=0.8] () at (15.5, 0.7166666666666668) {$n$};
\draw (14.0, 5.583333333333333) .. controls (14.0, 1.4166666666666665) .. (14.0, 1.4166666666666665);
\draw (16.0, 3.2500000000000004) .. controls (16.0, 1.4166666666666665) .. (16.0, 1.4166666666666665);
\draw (15.0, 0.9166666666666667) .. controls (15.0, 0.0) .. (15.0, 0.0);
\draw (13.5, 5.583333333333334) -- (14.5, 5.583333333333334) -- (14.5, 6.083333333333334) -- (13.5, 6.083333333333334) -- (13.5, 5.583333333333334);
\node [scale=0.8] () at (14.0, 5.833333333333334) {one};
\draw (15.5, 3.25) -- (16.5, 3.25) -- (16.5, 3.75) -- (15.5, 3.75) -- (15.5, 3.25);
\node [scale=0.8] () at (16.0, 3.5) {two};
\draw (13.5, 0.9166666666666667) -- (16.5, 0.9166666666666667) -- (16.5, 1.4166666666666667) -- (13.5, 1.4166666666666667) -- (13.5, 0.9166666666666667);
\node [scale=0.8] () at (15.0, 1.1666666666666667) {plus};
\end{tikzpicture} $\quad \mapsto \quad 3$
\end{center}
where the first step is snake removal and the second is function evaluation as in appendix~\ref{3-cartesian}.
\end{example}

\section{tensor.py}\label{5-tensor}

Let $\mathbf{Mat}_\bb{S}$ be the category with objects the natural numbers and arrows $m \to n$ the matrices $[n] \times [m] \to \bb{S}$ for a commutative semiring $\bb{S}$, with Kronecker product as tensor.
$\mathbf{Mat}_\bb{S}$ is compact closed, i.e. both symmetric monoidal and rigid. It is furthermore self-dual, i.e. objects are isomorphic to their adjoints.
For $\bb{S} = \bb{B}$, we get a category equivalent to finite sets and relations with Cartesian product.
For $\bb{S} = \bb{C}$, it is equivalent to finite-dimensional complex vector spaces and linear maps with the usual tensor product.
In practice, it is more convenient to consider an equivalent category $\mathbf{Tensor}_\bb{S}$ where objects are lists of natural numbers and adjoints are given by list reversal.
Arrows $(m_1, \dots, m_k) \to (n_1, \dots, n_{k'})$ are tensors of order $k + k'$, i.e. matrices $m_1 \times \dots \times m_k \to n_1 \times \dots \times n_{k'}$.

\begin{class}\normalfont\texttt{Dim(n\_0, ..., n\_k)}
is a subclass of \py{rigid.Ty} generated by natural numbers.
\end{class}

\begin{class}\normalfont\texttt{Tensor(dom, cod, array)}
is a subclass of \py{rigid.Box} given by \py{Dim}-instances \py{dom} and \py{cod} and a \texttt{numpy} \cite{VanderWaltEtAl11} \py{array} of the appropriate shape.
\py{then} and \py{tensor} are both implemented using \py{numpy.tensordot}, \py{cups} and \py{caps} return a reshaped identity.
\end{class}

\begin{class}\normalfont\texttt{TensorFunctor(ob, ar)}
is a subclass of \py{rigid.Functor} where \py{ob} and \py{ar} are mappings from \py{Ty} to \py{Dim} and from \py{Box} to \py{Tensor} respectively.
\end{class}

\begin{remark}
All the methods of the \py{Tensor} class are writen in \normalfont{\texttt{jax.numpy}}, the subset of Python+\normalfont{\texttt{numpy}} that supports automatic differentiation with \normalfont{\texttt{jax}} \cite{Google/jax20}.
\end{remark}

\begin{example}
Tensor networks can be defined as diagrams with a functor into tensors, contraction is given by functor application. They have been applied to both condensed matter physics and machine learning, see \cite{Orus14} for an introduction.
Interfacing \normalfont{DisCoPy} with tensor network tools such as \cite{KossaifiEtAl18,HauschildPollmann18,RobertsEtAl19} is left for future work.
\end{example}

\begin{example}
Relational databases can be defined as Boolean tensors: a table with $k$ columns is a state $1 \to (n_1, \dots, n_k)$ in $\mathbf{Tensor}_\bb{B}$.
Conjunctive queries are diagrams, where query containment gives the structure of a free Cartesian bicategory, see \cite{BonchiEtAl18}.
Query evaluation over a relational database is the application a functor into Boolean tensors.
\end{example}

\begin{example}\label{example-discocat}
The distributional compositional (DisCo) models of Coecke et al.
\cite{ClarkEtAl08,ClarkEtAl10} can be defined as functors
$F : \mathbf{G} \to \mathbf{Tensor}_\bb{S}$ from the rigid category $\mathbf{G}$ generated by
a pregroup grammar (see example~\ref{example-pregroup}) into tensors, i.e. they map pregroup types $t \in \mathbf{G}$ to dimensions $F(t) \in \bb{N}^\star$ and dictionnary entries $w \to t$ to tensors of shape $F(t)$.
When $F(s) = 1$, the meaning for a grammatical sentence $g : w_1 \dots w_n \to s$ is a scalar $F(g) \in \bb{S}$ which can be computed as the contraction of a tensor network.
DisCo models into real vector spaces, i.e. with $\bb{S} = \bb{R}$, received experimental support, see \cite{GrefenstetteSadrzadeh11,KartsaklisEtAl12,KartsaklisEtAl13}.
Relational DisCo models, i.e. with $\bb{S} = \bb{B}$, have been applied to question answering, see \cite{CoeckeEtAl18a, DeFeliceEtAl19a}.
\end{example}

\section{circuit.py}\label{6-circuit}

Quantum circuits are a standard model for quantum computation.
They form the arrows of a PROP, i.e. a symmetric monoidal category generated by one object, called a qubit.
We define $\mathbf{Circ}$ as the free PROP generated by
$n$-qubit gates $g : n \to n$,
scalars $\{s : 0 \to 0\}_{s \in \bb{C}}$,
post-selection $\{\mathtt{bra}_i : 1 \to 0\}_{i \in \{0, 1\}}$
and preparation $\{\mathtt{ket}_i : 0 \to 1\}_{i \in \{0, 1\}}$ of ancilla qubits in the computational basis.
Circuit evaluation is defined as a monoidal functor $\mathtt{eval} : \mathbf{Circ} \to \mathbf{Tensor}_\bb{C}$ which sends each gate $g : n \to n$ to its unitary matrix $\mathtt{eval}(g) : 2^n \to 2^n$.

Given the circuit for an $n$-qubit state $c : 0 \to n$, measurement results are a tensor $\mathtt{measure}(c) : 1 \to 2^n$ of non-negative reals in $\mathbf{Tensor}_{\bb{R}^+}$, computed using the Born rule.
Note that if a circuit contains scalars or post selection, the measurement results need not be a normalised probability distribution.

The quotient $\mathbf{Circ}_\sim$, where $c = c'$ iff $\mathtt{eval}(c) = \mathtt{eval}(c')$, is a compact-closed category.
Cups and caps are given by the (unnormalised) Bell effect and state, the snake equation implies the correctness of the teleportation protocol.
See \cite{CoeckeKissinger17} for an introduction to diagrammatic reasoning and quantum processes.

\begin{class}\normalfont\texttt{Circuit(dom, cod, boxes, offsets)}
is a subclass of \py{rigid.Diagram} with \py{PRO}-instances as \py{dom} and \py{cod}. It has methods \py{eval}, implemented as a \py{TensorFunctor}, and \py{measure} which computes the Born rule.
\end{class}

\begin{class}\normalfont\texttt{Bra(b\_0, ..., b\_n)} and \normalfont\texttt{Ket(b\_0, ..., b\_n)} are
subclasses of \py{Circuit} and \py{rigid.Box} given by a bitstring \py{b\_0, ..., b\_n}.
\end{class}

\begin{class}\normalfont\texttt{Gate(name, n\_qubits, array)}
is a subclass of \py{Circuit} and \py{rigid.Box} with instances \py{H}, \py{CX}, \py{SWAP}, etc. Phases are implemented as subclasses \py{Rx} and \py{Rz}.
\end{class}

\begin{class}\normalfont\texttt{CircuitFunctor(ob, ar)}
is a subclass of \py{rigid.Functor} where \py{ob} and \py{ar} are mappings from \py{Ty} to \py{PRO} and from \py{Box} to \py{Circuit} respectively.
\end{class}

The methods \py{to\_tk} and \py{from\_tk} translate back and forth between DisCoPy's \py{Circuit} class and that of t$\vert$ket$\rangle$ \cite{SivarajahEtAl20}, which can then be compiled and executed on quantum hardware or simplified using \texttt{pyzx} \cite{KissingervandeWetering19}.
Note that in the translation from DisCoPy diagrams to the directed acyclic graphs of t$\vert$ket$\rangle$, we treat the \py{SWAP} gate as a logical gate, i.e. it simply renames the two qubits.
In the other direction, we introduce \py{SWAP} gates whenever a t$\vert$ket$\rangle$ gate is applied to non-adjacent qubits.
Thus, \py{from\_tk(c.to\_tk())} is equal to the original circuit \py{c} up to the axioms of symmetric monoidal categories.

\begin{example}
The quantum algorithms for natural language processing (NLP) of \cite{ZengCoecke16} can be defined as rigid functors $\mathbf{G} \to \mathbf{Circ}$ from a pregroup grammar (see examples~\ref{example-pregroup} and \ref{example-discocat}) to the category of circuits.
See \cite{Coecke19} for a discussion of distributional compositional models for NLP on quantum hardware.
A proof-of-concept was implemented using \normalfont{DisCoPy}, see the notebook of the first experiments \href{https://github.com/oxford-quantum-group/discopy/blob/master/notebooks/qnlp-experiment.ipynb}{here} and \href{https://medium.com/cambridge-quantum-computing/quantum-natural-language-processing-748d6f27b31d}{there} \cite{Meichanetzidis20} for more details.
\end{example}

{\footnotesize
\bibliographystyle{eptcs}
\bibliography{discocat-complexity}}

\appendix
\section{cartesian.py}\label{3-cartesian}

This appendix describes \py{cartesian.Diagram} and \py{PythonFunctor}, an
implementation of \emph{Lawvere theories} and their models.
We first give a short introduction to symmetric monoidal categories (SMC), PROPs and functorial semantics.

A (strict) monoidal category $\mathbf{C}$ is \emph{symmetric} when it comes equipped
with a natural transformation $\sigma_{x, y} : x \otimes y \to y \otimes x$
such that $\sigma_{y, x} \circ \sigma_{x, y} = 1_{x \otimes y}$ (involution)
and $\sigma_{x, y \otimes z} =
(1_{z} \otimes \sigma_{x, z}) \circ (\sigma_{x, y} \otimes 1_{z})$ (hexagon).
A component $\sigma_{x, y}$ is depicted as a swap of the wires for $x$ and $y$.
The free symmetric monoidal category $\mathbf{SMC}(\Sigma)$ generated by a monoidal
signature $\Sigma$ can be defined as a quotient $\mathbf{MC}(\Sigma') / \cal{S}$ of
the free monoidal category generated by the disjoint union
$\Sigma' = \Sigma + \set{\sigma_{x, y}}_{x, y \in \Sigma_0}$. The relation
$\cal{S}$ is generated by the rules for involution and naturality:
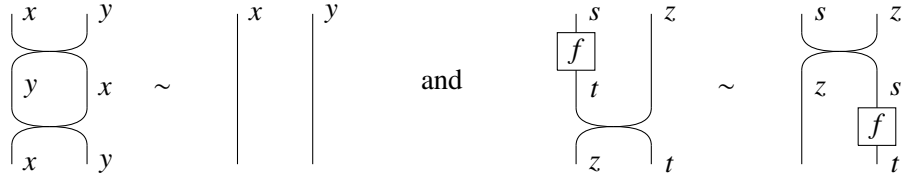
\begin{center}
\begin{tikzpicture}[baseline=(O.base)]
\node (O) at (0, 1.0) {};
\node () at (0.25, 2.0) {$x$};
\node () at (1.25, 2.0) {$y$};
\node () at (0.25, 1.0) {$y$};
\node () at (1.25, 1.0) {$x$};
\node () at (0.25, 0.0) {$x$};
\node () at (1.25, 0.0) {$y$};
\draw [out=-90, in=90] (0, 2.0) to (0, 1.75);
\draw [out=-90, in=90] (1, 2.0) to (1, 1.75);
\draw [out=180, in=90] (0.5, 1.5) to (0.0, 1.25);
\draw [out=0, in=90] (0.5, 1.5) to (1.0, 1.25);
\draw [out=-90, in=180] (0, 1.75) to (0.5, 1.5);
\draw [out=-90, in=0] (1, 1.75) to (0.5, 1.5);
\draw [out=-90, in=90] (0.0, 1.25) to (0.0, 0.75);
\draw [out=-90, in=90] (1.0, 1.25) to (1.0, 0.75);
\draw [out=180, in=90] (0.5, 0.5) to (0.0, 0.25);
\draw [out=0, in=90] (0.5, 0.5) to (1.0, 0.25);
\draw [out=-90, in=180] (0.0, 0.75) to (0.5, 0.5);
\draw [out=-90, in=0] (1.0, 0.75) to (0.5, 0.5);
\draw [out=-90, in=90] (0.0, 0.25) to (0.0, 0.0);
\draw [out=-90, in=90] (1.0, 0.25) to (1.0, 0.0);
\node () at (2, 1.0) {$\sim$};
\node () at (3.25, 2.0) {$x$};
\node () at (4.25, 2.0) {$y$};
\draw [out=-90, in=90] (3, 2.0) to (3, 0.0);
\draw [out=-90, in=90] (4, 2.0) to (4, 0.0);
\end{tikzpicture}
\qquad and \qquad \quad
\begin{tikzpicture}[baseline=(O.base)]
\node (O) at (0, 1.0) {};
\node () at (0.25, 2.0) {$s$};
\node () at (1.25, 2.0) {$z$};
\node () at (0.25, 1.0) {$t$};
\node () at (0.25, 0.0) {$z$};
\node () at (1.25, 0.0) {$t$};
\draw [out=-90, in=90] (0, 2.0) to (0, 1.75);
\draw [out=-90, in=90] (1, 2.0) to (1, 0.75);
\draw [out=-90, in=90] (0.0, 1.25) to (0.0, 0.75);
\draw [out=180, in=90] (0.5, 0.5) to (0.0, 0.25);
\draw [out=0, in=90] (0.5, 0.5) to (1.0, 0.25);
\draw [out=-90, in=180] (0.0, 0.75) to (0.5, 0.5);
\draw [out=-90, in=0] (1, 0.75) to (0.5, 0.5);
\draw [out=-90, in=90] (0.0, 0.25) to (0.0, 0.0);
\draw [out=-90, in=90] (1.0, 0.25) to (1.0, 0.0);
\draw (-0.25, 1.25) -- (0.25, 1.25) -- (0.25, 1.75) -- (-0.25, 1.75) -- (-0.25, 1.25);
\node () at (0.0, 1.5) {$f$};
\node () at (2, 1.0) {$\sim$};
\node () at (3.25, 2.0) {$s$};
\node () at (4.25, 2.0) {$z$};
\node () at (3.25, 1.0) {$z$};
\node () at (4.25, 1.0) {$s$};
\node () at (4.25, 0.0) {$t$};
\draw [out=-90, in=90] (3, 2.0) to (3, 1.75);
\draw [out=-90, in=90] (4, 2.0) to (4, 1.75);
\draw [out=180, in=90] (3.5, 1.5) to (3.0, 1.25);
\draw [out=0, in=90] (3.5, 1.5) to (4.0, 1.25);
\draw [out=-90, in=180] (3, 1.75) to (3.5, 1.5);
\draw [out=-90, in=0] (4, 1.75) to (3.5, 1.5);
\draw [out=-90, in=90] (3.0, 1.25) to (3.0, 0.0);
\draw [out=-90, in=90] (4.0, 1.25) to (4.0, 0.75);
\draw [out=-90, in=90] (4.0, 0.25) to (4.0, 0.0);
\draw (3.75, 0.25) -- (4.25, 0.25) -- (4.25, 0.75) -- (3.75, 0.75) -- (3.75, 0.25);
\node () at (4.0, 0.5) {$f$};
\end{tikzpicture}
\end{center}
for all $x, y, z \in \Sigma_0$ and $f : s \to t$ in $\Sigma'$. Note that
$s, t \in \Sigma_0^\star$ may be of arbitrary length, in which case $\sigma_{s, z}$ and $\sigma_{t, z}$ are defined as ladders of swaps,
i.e. the symmetry for compound types $\sigma_{x \otimes y, z}$ and
$\sigma_{x, y \otimes z}$ is defined inductively by:
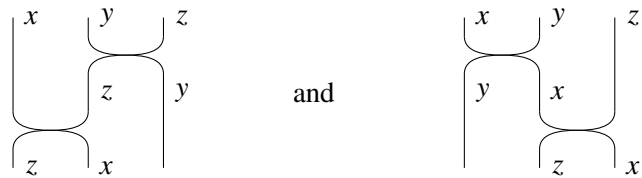
\begin{center}
\begin{tikzpicture}[baseline=(O.base)]
\node (O) at (0, 1.0) {};
\node () at (0.25, 2.0) {$x$};
\node () at (1.25, 2.0) {$y$};
\node () at (2.25, 2.0) {$z$};
\node () at (1.25, 1.0) {$z$};
\node () at (2.25, 1.0) {$y$};
\node () at (0.25, 0.0) {$z$};
\node () at (1.25, 0.0) {$x$};
\draw [out=-90, in=90] (0, 2.0) to (0, 0.75);
\draw [out=-90, in=90] (1, 2.0) to (1, 1.75);
\draw [out=-90, in=90] (2, 2.0) to (2, 1.75);
\draw [out=180, in=90] (1.5, 1.5) to (1.0, 1.25);
\draw [out=0, in=90] (1.5, 1.5) to (2.0, 1.25);
\draw [out=-90, in=180] (1, 1.75) to (1.5, 1.5);
\draw [out=-90, in=0] (2, 1.75) to (1.5, 1.5);
\draw [out=-90, in=90] (1.0, 1.25) to (1.0, 0.75);
\draw [out=-90, in=90] (2.0, 1.25) to (2.0, 0.0);
\draw [out=180, in=90] (0.5, 0.5) to (0.0, 0.25);
\draw [out=0, in=90] (0.5, 0.5) to (1.0, 0.25);
\draw [out=-90, in=180] (0, 0.75) to (0.5, 0.5);
\draw [out=-90, in=0] (1.0, 0.75) to (0.5, 0.5);
\draw [out=-90, in=90] (0.0, 0.25) to (0.0, 0.0);
\draw [out=-90, in=90] (1.0, 0.25) to (1.0, 0.0);
\node () at (4, 1.0) {and};
\node () at (6.25, 2.0) {$x$};
\node () at (7.25, 2.0) {$y$};
\node () at (8.25, 2.0) {$z$};
\node () at (6.25, 1.0) {$y$};
\node () at (7.25, 1.0) {$x$};
\node () at (7.25, 0.0) {$z$};
\node () at (8.25, 0.0) {$x$};
\draw [out=-90, in=90] (6, 2.0) to (6, 1.75);
\draw [out=-90, in=90] (7, 2.0) to (7, 1.75);
\draw [out=-90, in=90] (8, 2.0) to (8, 0.75);
\draw [out=180, in=90] (6.5, 1.5) to (6.0, 1.25);
\draw [out=0, in=90] (6.5, 1.5) to (7.0, 1.25);
\draw [out=-90, in=180] (6, 1.75) to (6.5, 1.5);
\draw [out=-90, in=0] (7, 1.75) to (6.5, 1.5);
\draw [out=-90, in=90] (6.0, 1.25) to (6.0, 0.0);
\draw [out=-90, in=90] (7.0, 1.25) to (7.0, 0.75);
\draw [out=180, in=90] (7.5, 0.5) to (7.0, 0.25);
\draw [out=0, in=90] (7.5, 0.5) to (8.0, 0.25);
\draw [out=-90, in=180] (7.0, 0.75) to (7.5, 0.5);
\draw [out=-90, in=0] (8, 0.75) to (7.5, 0.5);
\draw [out=-90, in=90] (7.0, 0.25) to (7.0, 0.0);
\draw [out=-90, in=90] (8.0, 0.25) to (8.0, 0.0);
\end{tikzpicture}
\end{center}
The symmetry for the empty type is defined to be the identity
$\sigma_{x, 1} = \sigma_{1, x} = 1_x$.
An SMC where the tensor is the Cartesian product and the unit is terminal
(i.e. a category with finite products) is called a \emph{Cartesian category}.
Equivalently, an SMC is Cartesian when objects carry a natural
commutative comonoid structure \cite[6.1]{Selinger10}.
Given a monoidal signature $\Sigma$, the free Cartesian category
$\mathbf{CC}(\Sigma)$ is the quotient
$\mathbf{MC}(\Sigma'') / (\cal{S} + \cal{P})$ for
$\Sigma'' = \Sigma' + \set{\mu_{x} : x \to x \otimes x}_{x \in \Sigma_0} + \set{\epsilon_x : x \to 1}_{x \in \Sigma_0}$.
The components $\mu_{x}$ and $\epsilon_{x}$ are depicted as wire splitting and
ending respectively. The comonoids of non-atomic
types inductively. That for the unit is the identity and the comonoid
of $x \otimes y$ is given by:
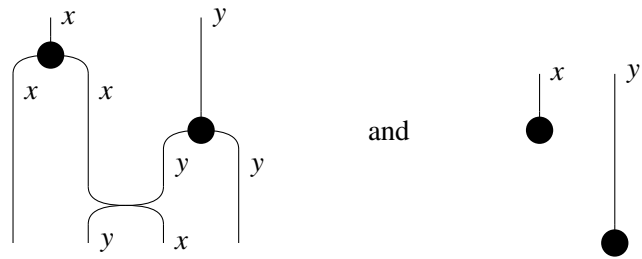
\begin{center}
\begin{tikzpicture}[baseline=(O.base)]
\node (O) at (0, 1.5) {};
\node () at (0.75, 3.0) {$x$};
\node () at (2.75, 3.0) {$y$};
\node () at (0.25, 2.0) {$x$};
\node () at (1.25, 2.0) {$x$};
\node () at (2.25, 1.0) {$y$};
\node () at (3.25, 1.0) {$y$};
\node () at (1.25, 0.0) {$y$};
\node () at (2.25, 0.0) {$x$};
\draw [out=-90, in=90] (0.5, 3.0) to (0.5, 2.75);
\draw [out=-90, in=90] (2.5, 3.0) to (2.5, 1.75);
\draw [out=180, in=90] (0.5, 2.5) to (0.0, 2.25);
\draw [out=0, in=90] (0.5, 2.5) to (1.0, 2.25);
\draw [out=-90, in=90] (0.5, 2.75) to (0.5, 2.5);
\draw [out=-90, in=90] (0.0, 2.25) to (0.0, 0.0);
\draw [out=-90, in=90] (1.0, 2.25) to (1.0, 0.75);
\draw [out=180, in=90] (2.5, 1.5) to (2.0, 1.25);
\draw [out=0, in=90] (2.5, 1.5) to (3.0, 1.25);
\draw [out=-90, in=90] (2.5, 1.75) to (2.5, 1.5);
\draw [out=-90, in=90] (2.0, 1.25) to (2.0, 0.75);
\draw [out=-90, in=90] (3.0, 1.25) to (3.0, 0.0);
\draw [out=180, in=90] (1.5, 0.5) to (1.0, 0.25);
\draw [out=0, in=90] (1.5, 0.5) to (2.0, 0.25);
\draw [out=-90, in=180] (1.0, 0.75) to (1.5, 0.5);
\draw [out=-90, in=0] (2.0, 0.75) to (1.5, 0.5);
\draw [out=-90, in=90] (1.0, 0.25) to (1.0, 0.0);
\draw [out=-90, in=90] (2.0, 0.25) to (2.0, 0.0);
\node [circle, fill=black] () at (0.5, 2.5) {};
\node [circle, fill=black] () at (2.5, 1.5) {};
\node () at (5.0, 1.5) {and};
\node () at (7.25, 2.25) {$x$};
\node () at (8.25, 2.25) {$y$};
\draw [out=-90, in=90] (7.0, 2.25) to (7.0, 1.75);
\draw [out=-90, in=90] (8.0, 2.25) to (8.0, 0.25);
\draw [out=-90, in=90] (7.0, 1.75) to (7.0, 1.5);
\draw [out=-90, in=90] (8.0, 0.25) to (8.0, 0.0);
\node [circle, fill=black] () at (7.0, 1.5) {};
\node [circle, fill=black] () at (8.0, 0.0) {};
\end{tikzpicture}
\end{center}
The relation $(\cal{S} + \cal{P})$ is given by the axioms for commutative
comonoids plus naturality of symmetry, coproduct and counit for each generating
arrow.
A \emph{PROP} (PROduct and Permutation) \cite{Lack04} is an SMC generated by one object, a \emph{Lawvere theory} \cite{Lawvere63} is a Cartesian category generated one object.

Lawvere theories are implemented by the class \py{cartesian.Diagram}, a subclass of\linebreak \py{monoidal.Diagram} with static methods \py{swap}, \py{copy} and \py{delete} implementing the structural morphisms.
\py{cartesian.Box} is a subclass of \py{monoidal.Box} where each instance holds a Python function with natural numbers as domain and codomain. A Cartesian box \py{f} with \py{f.dom, f.cod == (m, n)} sends $m$-tuples to $n$-tuples, it can be defined in the standard syntax for Python functions using the decorator \py{@disco(m, n)}.
\py{Swap}, \py{Copy} and \py{Del} are subclasses of \py{cartesian.Box} which implement the symmetry and comonoid on the generating object.
\py{cartesian.Functor} is a subclass of \py{monoidal.Functor} which preserves symmetry and product.
\py{Function} is a subclass of \py{cartesian.Box} where \py{then} and \py{tensor} are overriden by function composition and tuple concatenation.
The \py{PythonFunctor} class implements Cartesian functors into \py{Function}, i.e. it maps the formal composition of diagrams to the concrete composition of functions.
Note that when \py{f} and \py{g} have side-effects, the tensor \py{f @ g} is in
general different from \py{Id(f.dom) @ g >> f @ Id(g.cod)}. Thus, \py{Function}
is closer to the implementation of a premonoidal than a monoidal category,
see example~\ref{example-2}.

\begin{example}
The Lawvere theory $\mathbf{F}$ with no generating arrows is the
opposite of the category of finite sets with disjoint union as monoidal structure: diagrams $f : m \to n$ in
$\mathbf{F}$ correspond precisely to the graphs of the functions $f : [n] \to [m]$.
The subcategories of diagrams in $\mathbf{F}$ with no coproduct, no counit and no symmetry correspond to injective, surjective and monotone functions respectively.
The subcategory of $\mathbf{F}$ generated by symmetry alone, i.e. the free SMC generated by one object, corresponds to bijections.
The normal form for diagrams in $\mathbf{F}$ is given by decomposing any function into a surjection and a monotone injection, see \cite[Theorem~3]{Lafont03}.
\end{example}

\begin{example}
For a semiring $\bb{S}$, the category of matrices over $\bb{S}$ with direct sum as monoidal product can be defined as a Lawvere theory generated by a commutative monoid $+ : 2 \to 1$ with unit $e : 0 \to 1$ and scalars $\set{s : 1 \to 1}_{s \in \bb{S}}$.
The relations are given by the axioms for semirings, see \cite{WadsleyWoods15}.
The naturality relation for the comonoid with respect to $+$ and $e$ are called the bialgebra laws, they have applications in control and network theory \cite{BaezErbele14,BaezEtAl18}.
\end{example}

\begin{example}
Let $\mathbf{NN}$ be the Lawvere
theory generated by sum $+ : 2 \to 1$, activation $a : 1 \to 1$,
finite sets of weights $\set{w_i : 1 \to 1}_{i \in W}$ and biases
$\set{b_i : 0 \to 1}_{i \in B}$. The arrows of $\mathbf{NN}$ are the diagrams
of neural network architectures.
Given a set of parameters $\theta : W + B \to \bb{R}$, an implementation is
given by the model $F_\theta : \mathbf{NN} \to \mathbf{Set}$ such that
$F_\theta(1) = \bb{R}$,
$F_\theta(b_i) = \set{() \mapsto \theta(b_i)}$,
$F_\theta(w_i) = \set{x \mapsto \theta(w_i) \cdot x}$,
$F_\theta(+) = \set{(x, y) \mapsto x + y}$
and $F_\theta(a) : \bb{R} \to \bb{R}$ is a non-linearity such as sigmoid.
A loss $l : (\bb{R}^m \to \bb{R}^n) \to \bb{R}$ for an architecture $f : m \to n$ in $\mathbf{NN}$ takes an implementation and
returns a real number encoding its success at some data-driven task.
The gradient of $\theta \mapsto l(F_\theta(f))$ can be computed by
back-propagation over the architecture $f$, using automatic differentiation tools such as $\mathtt{jax}$ \cite{Google/jax20}. The back-propagation algorithm is
itself part of a functor $\mathbf{NN} \to \mathbf{Learn}$, where $\mathbf{Learn}$ is a
Cartesian category of supervised learning algorithms, see \cite{FongEtAl17, FongJohnson19}.
\end{example}

\section{drawing.py}\label{a-draw}

We reformulate some of the definitions and results from Joyal, Street \cite{JoyalStreet88}.

\begin{definition}
A \emph{topological graph}, also called 1d cell complex, is a tuple $(\Gamma, \Gamma_0, \Gamma_1)$ of a Hausdorff space $\Gamma$ and a closed, discrete subset $\Gamma_0 \sub \Gamma$ of nodes such that $\Gamma - \Gamma_0 = \bigsqcup \Gamma_1$ is the sum of its connected components $\Gamma_1$, called edges, each homeomorphic to an open intervals with boundary in $\Gamma_0$.
\end{definition}

\begin{definition}
A planar graph between two real numbers $a < b$ is a finite topological graph $\Gamma$ embedded in $\bb{R} \times [a, b]$ such that every $x \in \Gamma \cap (\bb{R} \times \set{a, b})$ is a node in $\Gamma_0$ and belongs to the closure of exactly one edge in $\Gamma_1$.
\end{definition}

The points in $\Gamma_0 \cap (\bb{R} \times \set{a, b})$ are called outer nodes,
they are the boundary of the domain and codomain $\mathtt{dom}(\Gamma), \mathtt{cod}(\Gamma) \in \Gamma_1^\star$ of the planar graph $\Gamma$.
The points $f \in \Gamma_0 - (\bb{R} \times \set{a, b})$ are called inner nodes, they have a domain and codomain $\mathtt{dom}(f), \mathtt{cod}(f) \in \Gamma_1^\star$ given by the edges that have $f$ as boundary.
The composition $\Gamma \circ \Gamma'$ and tensor $\Gamma \otimes \Gamma'$ are defined by rescaling and pasting the two planar graphs $\Gamma, \Gamma'$ vertically and horizontally respectively, see \cite[§4]{JoyalStreet88}.

A planar graph $\Gamma$ is \emph{progressive}, or recumbent, when the second projection $e \to [a, b]$ is injective for every edge $e \in \Gamma_1$. Progressive planar graphs have no cups or caps.
A progressive planar graph is \emph{generic} when the projection $\Gamma_0 - (\bb{R} \times \set{a, b}) \to (a, b)$ is injective, i.e. there are no two inner nodes at the same height.
A deformation of planar graphs is a continuous map $h : \Gamma \times [0, 1] \to [a, b] \times \bb{R}$ such that 1) for all $t \in [0, 1]$, $h(-, t)$ is an embedding whose image is a planar graph, and 2) for all $x \in \Gamma_0$, if $h(x, t)$ is inner for some $t$ then it is inner for all values of $t \in [0, 1]$.
A deformation of planar graphs is progressive (generic) when for all $t \in [0, 1]$, $h(-, t)$ is progressive (generic).

The induced equivalence relation --- with $\Gamma_0 \sim \Gamma_1$ if and only if there is a deformation $h$ with $h(-, 0) = \Gamma_0$ and $h(-, 1) = \Gamma_1$ --- is a congruence with respect to composition and tensor of planar graphs,
i.e. if $\Gamma_0 \sim \Gamma_1$ and $\Omega_0 \sim \Omega_1$ then $\Gamma_0 \otimes \Omega_0 \sim \Gamma_1 \otimes \Omega_1$ and $\Gamma_0 \circ \Omega_0 \sim \Gamma_1 \circ \Omega_1$.
Furthermore, tensor and composition are associative and they respect the interchange law up to progessive deformation.
Thus, progessive planar graphs up to progessive deformation form a strict monoidal category \cite[Proposition~4]{JoyalStreet88}.

Given a monoidal signature $\Sigma$, a \emph{valuation} $v$ of a planar graph $\Gamma$ is a pair of functions $v_0 : \Gamma_1 \to \Sigma_0$ and $v_1 : \Gamma_0 - (\bb{R} \times \set{a, b}) \to \Sigma_1$ that send edges to objects and inner nodes to arrows, which commute with domain and codomain.
Progressive planar graphs valued in $\Sigma$ (up to progressive deformation) are the free monoidal category \cite[Theorem~5]{JoyalStreet88}.
We conjecture that generic planar graphs up to generic deformation are the free premonoidal category.

From the universal property, we know that the category of progressive planar graph is equivalent to the combinatorial definition of free monoidal categories given in section~\ref{2-monoidal}.
Concretely, this equivalence is witnessed by the following pair of algorithms translating between planar graphs and diagrams.

\vspace{-5pt}
\begin{algorithm}
\SetKwInOut{Input}{Input}\SetKwInOut{Output}{Output}
\Input{A progressive planar graph $\Gamma$}
\Output{\py{Diagram(dom, cod, boxes, offsets)}}
\BlankLine
\DontPrintSemicolon
Compute the connected components $\Gamma_1$ and their boundaries $\Gamma_0$.\;
Order the nodes by height and partition them into $\Gamma_0 = \py{dom} + \py{boxes} + \py{cod}$.\;
Find the start and end points $\Gamma_1 \to (\py{dom} + \py{boxes}) \times (\py{boxes} + \py{cod})$ for each edge, the preimage of this map gives the domain and codomain for each box.\;
Compute $\py{offsets}$ as the number of edges to the left of each box.\;
\caption{read}
\end{algorithm}
\vspace{-22pt}
\begin{algorithm}
\SetKwInOut{Input}{Input}\SetKwInOut{Output}{Output}
\Input{\py{Diagram(dom, cod, boxes, offsets)}}
\Output{A progressive planar graph $\Gamma$}
\BlankLine
\DontPrintSemicolon
Set $\Gamma_0 := \set{(i, 0) \ \vert \ i < \py{len(dom)}}$ and $\Gamma_1 := \emptyset$.\;
\For{\py{height, (box, offset) in enumerate(zip(boxes, offsets))}}{
Deform $\Gamma$ so that there is at least \py{len(box.cod) + 1} horizontal space between the edges of $\mathtt{cod}(\Gamma)$ at index \py{offset} and \py{offset + len(box.dom)}.\label{step-3}\;
Set $\Gamma_0 := \Gamma_0 + \{\py{box}\}$ for $\py{box} = (m, \py{height} +\frac{1}{2})$ with $m$ computed at step~\ref{step-3} and $\Gamma_1 := \Gamma_1 + \{x \to \py{box} \ \vert \ x \in
\py{box.dom}\} + \{\py{box} \to x \ \vert \ x \in
\py{box.cod}\}$.
}
Return $\Gamma := \Gamma + \{(i, j) \to (i, \py{max(len(boxes), 1)}) \ \vert \ (i, j) \in \mathtt{cod}(\Gamma)\}$.
\caption{draw}
\end{algorithm}

\end{document}